\documentstyle[12pt]{article}

  \newcommand{\const}{\rm const}
  \newcommand{\Var}{\rm Var}

\begin{document}

   \begin{center}

{\bf Tail estimations for normed sums of } \\

\vspace{4mm}

{\bf centered exchangeable random variables. }\\

\vspace{5mm}

{\bf M.R.Formica, \ E.Ostrovsky, \ L.Sirota. }

 \end{center}

 \vspace{5mm}

\ Universit\`{a} degli Studi di Napoli Parthenope, via Generale Parisi 13, Palazzo Pacanowsky, 80132,
Napoli, Italy. \\
e-mail: \ mara.formica@uniparthenope.it \\

 \ Department of Mathematics and Statistics, Bar-Ilan University,\\
59200, Ramat Gan, Israel. \\
e-mail: \ eugostrovsky@list.ru\\

 \ Department of Mathematics and Statistics, Bar-Ilan University,\\
59200, Ramat Gan, Israel. \\
e-mail: \ sirota3@bezeqint.net \\

\vspace{5mm}

\begin{center}

{\bf Abstract.}

\end{center}

\vspace{4mm}

 \ We derive in this short report  the exponential as well as power decreasing tail estimations for the sums of centered
exchangeable random variables, alike ones for the sums of the centered independent ones. \par

\vspace{5mm}

\par

\vspace{5mm}

 \section{Statement of problem. Notations.}

\vspace{5mm}

 \hspace{3mm}  Let $ \  (\Omega =\{\omega\}, \cal{B}, {\bf P} ) \ $ be certain non - trivial probability space with expectation $ \ {\bf E} \ $ and
variation $ \ {\bf \Var}. \ $   Let also $ \  (Y, \cal{N} )  \ $ be another measurable space.  \par

\vspace{4mm}

 \ {\bf Definition 1.1,} see \cite{Aldous}, \cite{Finetti 1},  \cite{Finetti 2},  \cite{Hewitt},   \cite{Lauritzen 1}, \cite{Lauritzen 2} at all. \par

 \vspace{3mm}

 \ The sequence of r.v.  $ \ \{\xi_i\}, \ i = 1,2,\ldots,n; \ n \le \infty \ $  with values in $ \ Y \ $ is said to be {\it exchangeable,} if for all the
permutations $ \ \pi \ = \ (\pi(1),\pi(2), \ \ldots,\pi(n)) \ $ of the set of numbers $ \  N = (1,2,\ldots,n),  \  $
 where $ \ n < \infty, \ $ the distribution of the random vector  $ \ (\xi_{\pi(1)}, \xi_{\pi(2)}, \ldots, \xi_{\pi(n)}) \ $ coincides
with the distribution of the source  random vector $ \ (\xi_{1}, \xi_{2}, \ldots, \xi_{n}): \ $

\vspace{3mm}

\begin{equation} \label{permut}
{\bf P} \left(\xi_{\pi(1)} \in A_1, \xi_{\pi(2)} \in A_2, \ldots, \xi_{\pi(n)} \in A_n \right) =
{\bf P} \left(\xi_1 \in A_1,\xi_{2} \in A_2, \ldots, \xi_{n} \in A_n \right)
\end{equation}
for all the tuples of the measurable sets $ \ A_j \in \cal{N}. \ $ The infinity sequence of the r.v. $ \ \{\xi_i\}, \ i \in 1,2, 3, \ldots \ $
is named exchangeable, if the relation (\ref{permut}) is satisfied for all the finite values $ \ n. \ $ \par

\vspace{3mm}

 \ Many examples of these type of the random variables, aside from the classical case of independent identical distributed r.v.,
  are described in the mentioned works; some applications in physic and statistic may be found in \cite{Kingman}. It was obtained
in particular the classical limit theorems: LLN, CLT, LIL etc. for the exchangeable r.v., still with the remainder term evaluation.\par

\vspace{4mm}

 \ {\bf Our aim in this short report is the non - asymptotic  estimation of the tail distribution, power as well as exponential decreasing,
  for the classical normed sums of such a random variables}

\begin{equation} \label{tail statement}
R_n(t) =  R[S(n)](t) \stackrel{def}{=} {\bf P} \left( \ n^{-1/2}|\sum_{i=1}^n \xi_i| \ge t \ \right), \ t \ge 1;
\end{equation}
and consequently

\begin{equation} \label{uniform tail statement}
\overline{R}(t) = \sup_n R[S(n)](t) \stackrel{def}{=} \sup_n {\bf P} \left( \ n^{-1/2}|\sum_{i=1}^n \xi_i| \ge t \ \right), \ t \ge 1;
\end{equation}

$$
S(n) \stackrel{def}{=} n^{-1/2} \sum_{i=1}^n \xi_i;
$$
of course, when $ \ Y = R^1, \ {\bf E} \xi_i = 0; \ \Var(\xi_i) \in (0,\infty); \ $ in the spirit  of the classical estimations belonging to
Dharmadhikari - Jogdeo, Rosenthal, Nagaev, Pinelis, Schechtman at all. See also \ \cite{ErOs}, \cite{Ermakov}, \cite{formicagiovamjom2015},
\cite{Formica Ostrovsky Sirota weak dep}, \cite{Ostrovsky1}, chapters 1,2.\par

\vspace{5mm}

 \section{Auxiliary fact: de Finetti representation.}

\vspace{5mm}

 \hspace{3mm} Let  now  in addition $ \ Y \ $ be metrisable separable complete space and $ \ \cal{N} \ $ be Borelian sigma field. We
 will use the following important fact: theorem  (representation) belonging to de Finetti for the exchangeable r.v., see  \cite{Finetti 1}, \cite{Finetti 2};
 and the more general proposition  \cite{Aldous}, \cite{Hewitt}.\par

 \vspace{4mm}

 \ {\bf Proposition 2.1.} There exists a probability measure $ \ \mu \ $ defined on the Borelian sigma field $ \ \cal{N} \ $  such that

\vspace{3mm}

\begin{equation} \label{Q represent}
{\bf P} \left( \ \xi_1 \in A_1,\xi_{2} \in A_2, \ldots, \xi_{n} \in A_n \ \right) =
\int_{\cal{N}}  Q(A_1) Q(A_2) \ldots Q(A_n) \ \mu(dQ).
\end{equation}

\vspace{4mm}

 \ Here in turn  $ \  Q(\cdot)  \ $  is certain probability measure. \par

\vspace{4mm}

 \ On the other words, if we introduce the r.v. $ \ M \ $ having the distribution $ \ \mu, \ $ then

\begin{equation} \label{Q condit}
{\bf P} \left(\xi_1 \in A_1,\xi_{2} \in A_2, \ldots, \xi_{n} \in A_n/M = Q \right) = Q(A_1) Q(A_2) \ldots Q(A_n),
\end{equation}
or in turn equally the r.v. $ \ \{\xi_i\} \ $ are {\it independent and identical distributed} under arbitrary {\it fixed} condition  $ \ M = Q; \ $
and in this condition they have the distribution  $ \ Q(\cdot). \ $ \par
 \ Herewith the measure $ \ M \ $ may be expressed as  almost everywhere  existing  limit

\begin{equation} \label{M limit}
M(A) = \lim_{n \to \infty} n^{-1} \sum_{i=1}^{\infty}\chi_A(\xi_i),
\end{equation}
where $ \ \chi_A(\cdot) \ $ is an indicator function for the Borelian set $ \ A. \ $ \par

\vspace{3mm}

 \ The inversely assertion to the proposition 2.1.  is obviously also true. Indeed, if for some random vector $  \  \{ \xi_i \}, \ i = 2,3,\ldots,n  \ $
 there holds true the representation (\ref{Q represent}), this vector is exchangeable. \par

\vspace{3mm}

 \ We suppose in addition that  the r.v. $ \  \{\xi_j\}  \ $ are conditional centered:

\begin{equation} \label{condit centered}
\forall Q \in {\cal N} \ \Rightarrow  {\bf E} \xi_i/Q = 0,
\end{equation}
and have a finite a.e. non - zero conditional variance (a second conditional moment)

\begin{equation} \label{condit variance}
\forall Q \in {\cal N} \ \Rightarrow  {\bf \Var} (\xi_i)/Q \in (0, \infty).
\end{equation}

\vspace{3mm}

 \ Denote yet

\begin{equation} \label{cond tails}
T[Q, n](t) \stackrel{def}{=} {\bf P} \left( \ n^{-1/2} \sum_{j=1}^n \xi_j > t   \ / Q \ \right), \ Q \in \cal{N},
\end{equation}
conditional tail of distribution if the normed sums, and consequently its uniform tail estimate

\begin{equation} \label{ supr cond tails}
\overline{T}[Q](t) \stackrel{def}{=} \sup_n T[Q, n](t).
\end{equation}
 \hspace{3mm}  Of course, under formulated above restrictions $ \ \overline{T}[Q](t) \to 0, \ t \to \infty.  \ $ \par

 \vspace{4mm}

\vspace{5mm}

 \section{Main result.}

\vspace{5mm}

 \hspace{3mm} {\bf Proposition 3.1.} It follows immediately from the well known complete probability formulae  that

\begin{equation} \label{main final  formula}
R_n(t) = \int_{\cal{N}} T[Q, n](t) \ \mu(dQ),
\end{equation}
and consequently

\begin{equation} \label{uniform main formula}
\overline{R}(t) \le \int_{\cal{N}} \overline{T}[Q](t) \ \mu(dQ).
\end{equation}

\vspace{5mm}

\section{Examples.}

\vspace{5mm}

 \hspace{3mm} There are a huge number works devoted to the estimations of the probability $ \ T[Q,n](t), \ $
see e.g. \cite{Buldygin},  \cite{Capone1}, \cite{Capone2}, \cite{ErOs}, \cite{Ermakov}, \cite{KosOs}, \cite{KozOsSir2017},
\cite{Liflyand}, \cite{Ostrovsky1}, \cite{Pekshir Shiryaev} etc.  For instance, let $ \ \{ \eta_i \},  \ i = 1,2, \ldots,n \ $
be a sequence of independent  identical distributed (i; i.d.) non trivial centered $ \ {\bf E} \eta_i = 0 \ $ random variables (r.v.). Set
as before

$$
\Theta_n := n^{-1/2} \sum_{i=1}^n \eta_i,
$$
and suppose in addition

$$
\exists \ m = \const  > 0 \ \Rightarrow \  {\bf P} ( \ |\eta_i| \ge u) \le \exp (- u^m), \ u \ge 0.
$$
 \ Then

$$
\exists c = c(m) \in (0,\infty) \ \Rightarrow  \sup_n {\bf P} (|\Theta_n|> u) \le \exp \left( - c(m) \ u^{\min(m,2)} \ \right), \ u \ge 0,
$$
and this estimate is essentially non - improvable as $ \ u \to \infty. \ $ \par

\vspace{4mm}

 \hspace{3mm} {\bf Example 4.1.} \ So,  it is reasonable to assume that the set $ \ \{Q\} \ $ may be identified with positive semi - axes
  $ \  \{ Q \}  = (0,\infty)  \ $ and let

\begin{equation} \label{Tail condit}
T[Q,n](t) \le \exp \left( \ - c_1 \ t^{\alpha} \ Q^{\beta} \ \right), \ c_1, \alpha,\beta \in (0,\infty);
\end{equation}

\begin{equation} \label{mu condit}
\mu(A) =  c_2  \ \int_A \ Q^{\gamma} \ \exp \left( \ - c_3 \ Q^{\kappa} \ \right) \ d Q, \ \gamma > - 1, \ \kappa > 0, \ c_3 > 0;
\end{equation}
where  $ \  A \ $ is Borelian subset of the  positive half line $ \ R_1: \ A \subset R_1 \ $  and of course (norming condition)

$$
c_2 = c_2(\kappa, c_3, \gamma) = \frac{\kappa \ c_3^{(\gamma + 1)/\kappa}}{\Gamma((\gamma + 1)/\kappa)},
$$
$ \ \Gamma(\cdot) \ $ is ordinary Gamma function. \par

\vspace{3mm}

 \ We get to the following expression for the estimation of the value  $ \ \overline{R}(t) \ $ from (\ref{uniform main formula}) under formulated
assumptions

\begin{equation} \label{tail expression}
\overline{R}(t) \le  c_2 \ \int_0^{\infty}  Q^{\gamma} \ \exp \left( \ - c_1 \ t^{\alpha} \ Q^{\beta} - c_3 Q^{\kappa} \ \right) \ d Q.
\end{equation}

 \vspace{4mm}

\begin{center}

 {\sc We came to the need to estimate the following  interest auxiliary integral}  \par

\end{center}

\vspace{3mm}

\begin{equation} \label{aux int}
I[\theta,g](t) \stackrel{def}{=} \int_0^{\infty} x^{\theta - 1} \ \exp (- t x - g(x)) \ dx,
\end{equation}
as $ \ t \to \infty, \ $ where $ \ \theta = \const > 0, \ g = g(x), \ x \ge 0 \ $ is non -  negative measurable
 function  such that $ \ g(0) = g(0+) = 0. \ $ \par

  \vspace{4mm}

 \ {\bf Lemma 4.1.} We conclude under formulated conditions as $ \ t \to \infty \ $

 \vspace{4mm}

\begin{equation} \label{key integral}
\lim_{t \to \infty}  \left\{ \ t^{\theta} I[\theta,g](t) \ \right\} = \Gamma(\theta).
\end{equation}

\vspace{4mm}

\ {\bf Proof. Upper estimation.}  We will use the classical saddle - point method, see e.g.
\cite{Fedoryuk}, \cite{Wong}. \par
\ We have  by virtue of the non - negativity of the function $ \ g = g(x) \ $
for all the sufficient greatest values $ \ t; \ $ say, for  the values $ \ t \ge 1 \ $

\begin{equation} \label{upp est}
I[\theta,g](t) \le \int_0^{\infty} x^{\theta - 1} \ \exp (- t x ) \ dx = \Gamma(\theta) \ t^{-\theta}.
\end{equation}

 \vspace{3mm}

 \ {\bf Lower estimate.} Let again $ \ t \ge 1. \ $ Let also $ \  \Delta = \Delta(t)  \ $ be certain numerical valued positive function
such that $ \ \lim_{t \to \infty} \overline{\Delta}(t) = 0;   \ $ where $ \  \overline{\Delta}(t) := \sup_{s \ge t} \Delta(s).   \ $ \par
 \ We have as $ \ t \to \infty \ $

$$
I[\theta,g](t) \ge \int_0^{\Delta} \exp( \ - tx - g(x) \ ) \ dx \ge \exp(-\overline{\Delta}) \cdot \int_0^{\Delta} e^{-tx} \ x^{\theta - 1} \ dx =
$$

$$
t^{-\theta}  \ \exp(-\overline{\Delta}) \ \int_0^{t \ \Delta} y^{\theta - 1} \ e^{-y} \ dy.
$$
\ It remains to choose for example $ \  \Delta = \Delta(t) := t^{-1/2},  \ $  to make sure of the fairness of our lemma. \par

\vspace{4mm}

 \ {\bf Remark 4.1.} One can generalize the proposition of this lemma: consider the following integral

$$
J[\theta, g, L](t) = \int_0^{\infty}  \ x^{\theta - 1} \ e^{-tx - g(x)} \ L(x) \ dx,
$$
where  in addition $ \ L = L(x) \ $ is bounded non - negative continuous {\it slowly varying}  at origin function; then as $ \ t \to \infty \ $

$$
\lim_{t \to \infty} \left\{ \ t^{\theta} J[\theta, g, L](t) \ \right\} = \int_0^{\infty} \ y^{\theta - 1} \ \exp(-y) \ L(y) \ dy.
$$

\vspace{4mm}

 \ Let us return to the source tail estimation (\ref{tail expression}). As long as, by virtue of the proposition of Lemma 4.1,  we have
 as $ \ t \to \infty \ $

\begin{equation} \label{tail transform begin}
  c_2 \ \int_0^{\infty}  Q^{\gamma} \ \exp \left( \ - c_1 \ t^{\alpha} \ Q^{\beta} - c_3 Q^{\kappa} \ \right) \ d Q \sim
\end{equation}

\begin{equation} \label{tail transform begin}
c_1^{-(\gamma + 1)/\beta} \  c_2  \ \beta^{-1} \ \Gamma((\gamma+ 1)/\beta) \ t^{ - \alpha (\gamma + 1)/\beta},
\end{equation}
following

\begin{equation} \label{key tail estim}
\overline{R}(t) \le C_4(\alpha,\beta,\kappa) \ t^{ - \alpha (\gamma + 1)/\beta}, \ t \ge 1, \ -
\end{equation}
the power tail  decay. \par

\vspace{4mm}

 \hspace{3mm} {\bf Example 4.2.} \  We  assume here again  that $ \  \{ Q \}  = (0,\infty)  \ $ and let now

\begin{equation} \label{Tail condit second}
T[Q,n](t) \le \exp \left( \ - c_1 \ t^{\alpha} \ Q^{ -\beta} \ \right), \ c_1, \alpha,\beta \in (0,\infty);
\end{equation}
and as above

\begin{equation} \label{mu condit sec}
\mu(A) =  c_2  \ \int_A \ Q^{\gamma} \ \exp \left( \ - c_3 \ Q^{\kappa} \ \right) \ d Q, \ \gamma > - 1, \ \kappa > 0, \ c_3 > 0;
\end{equation}

$$
c_2 = c_2(\kappa, c_3, \gamma) = \frac{\kappa \ c_3^{(\gamma + 1)/\kappa}}{\Gamma((\gamma + 1)/\kappa)}.
$$

\vspace{3mm}

 \ We get to the following expression for the estimation of the  new value  $ \ \overline{R}(t) \ $ from (\ref{uniform main formula}) under formulated
in this example assumptions  $ \  \overline{R}(t) \le  R_0(t) = R_0[\alpha, \ \beta, \ \gamma, \ \kappa, \ c_1,\ c_2, \ c_3](t),    \ $  where

\begin{equation} \label{new tail expression}
R_0(t) \stackrel{def}{=}   c_2 \ \int_0^{\infty}  Q^{\gamma} \ \exp \left( \ - c_1 \ t^{\alpha} \ Q^{ -\beta} - c_3 Q^{\kappa} \ \right) \ d Q.
\end{equation}

\vspace{3mm}

 \ One can apply for the investigation of the last integral (\ref{new tail expression})  as $ \ t \to \infty \ $ the classical saddle - point method,
 to obtain the logarithmical exact estimation

\begin{equation} \label{log exact}
\overline{R}(t) \le c_6(\alpha, \beta,\kappa,\gamma) \ \exp \left\{ \ - c_7(\alpha, \beta,\kappa,\gamma) \ t^{\alpha \kappa/(\beta + \kappa)} \ \right\}, \ t \ge 1,
\end{equation}
for some finite non - zero "constants" $ \ c_6, c_7; \ $ the exponential tail decay. \par

\vspace{5mm}

 \  More detail computations by means of  the saddle - point (Laplace) method
 show us the following  evaluation as $ \ t \to \infty \ $

$$
R_0(t)  \sim  \sqrt{2 \ \pi} \ c_2 \ c_{10}  \ t^{c_{11}} \ (A + B)^{-1/2} \ \times
$$

$$
 \exp \ \left\{ \ - c_{10} \ \left(A^{-1} + B^{-1} \right) \ t^{\alpha \ \kappa/(\beta + \kappa)}  \ \right\},
$$
 where

$$
 A = \frac{\beta}{\gamma + 1}, \ \hspace{3mm} B = \frac{\kappa}{\gamma + 1},
$$

$$
c_{10} = (\gamma + 1)^{-1} \cdot \left[ \ c_1^{\kappa} \ c_3^{\beta} \ \beta^{\kappa} \ \kappa^{\beta} \ \right]^{1/(\beta + \kappa)}, \hspace{3mm}
 c_{11} = \alpha \cdot \left[ \ \frac{2 \gamma + 2 - \kappa}{2(\beta + \kappa)} \ \right] - 1.
$$

 \vspace{5mm}

 \section{Concluding remarks.}

 \vspace{5mm}

 \hspace{3mm} \ It is interest in our opinion to obtain the multivariate generalization of our results, as well as
 to consider another examples.\par

\vspace{6mm}

\vspace{0.5cm} \emph{Acknowledgement.} {\footnotesize The first
author has been partially supported by the Gruppo Nazionale per
l'Analisi Matematica, la Probabilit\`a e le loro Applicazioni
(GNAMPA) of the Istituto Nazionale di Alta Matematica (INdAM) and by
Universit\`a degli Studi di Napoli Parthenope through the project
\lq\lq sostegno alla Ricerca individuale\rq\rq .\par

\end{document}